\documentclass[12pt]{article}
\input{epsf.sty}
\usepackage{epsfig}
\usepackage{times}
\newcommand{\onetom}{1,2,\cdots,m}

\newtheorem{mythm}{Theorem}
\newtheorem{mycol}{Corollary}
\newtheorem{mylem}{Lemma}

\newtheorem{mydef}{Definition}
\newtheorem{myrem}{Remark}
\def\sqr#1#2{{\vcenter{\vbox{\hrule height.#2pt \hbox{\vrule width.
#2pt height#1pt \kern#1pt \vrule} \hrule height.#2pt}}}}

\begin{document}

\begin{center}
{\large\bf Network synchronization with an adaptive strength}
\\[0.2in]
\begin{center}
XIWEI LIU AND TIANPING CHEN

Lab. of Nonlinear Mathematics Science, Institute of Mathematics,
Fudan University, Shanghai, 200433, P.R.China.\\
\indent ~~Corresponding author: Tianping Chen.
Email:tchen@fudan.edu.cn
\end{center}

\end{center}

\begin{abstract}
In this paper, new schemes to synchronize linearly or nonlinearly
coupled chaotic systems with an adaptive coupling strength are
proposed. Unlike other adaptive schemes, which synchronize coupled
chaotic systems to a special trajectory ( or an equilibrium point)
of the uncoupled node by adding negative feedbacks adaptively; here,
adaptive schemes for the coupling strength are used to synchronize
coupled chaotic systems without knowing the synchronization
trajectory. Moreover, in many applications, the state variables are
not observable; instead, some functions of the states can be
observed. How to synchronize coupled systems with the observed data
is of great significance. In this paper, synchronization of
nonlinearly coupled chaotic systems with an adaptive coupling
strength is also discussed. The validity of those schemes are proved
rigorously. Moreover, simulations show that by choosing some proper
parameter $\alpha$, the coupling strength obtained by adaptation
could be much smaller. It means that chaos oscillators can be easily
synchronized with a very weak coupling.


\end{abstract}

\pagestyle{plain} \pagenumbering{arabic}

\section{Introduction}\quad
Recently, an increasing interest has been devoted to the study of
complex networks (see Strogaze 2001 - Wang et.al 2002). Among
them, synchronization of coupled complex networks has received
more attentions, because synchronization not only can explain many
natural phenomena (see Mirollo and Strogatz 1990), but also have
many applications, such as image processing, secure communication
(see Wei and Jia 2002 - Lu and Chen 2004) and so on.

Generally, linearly coupled systems can be described as:
\begin{eqnarray}\label{synchro}
\dot{x}_i(t)=f(x_i(t),t)+c\sum\limits_{j\not =i}a_{ij}\Gamma
\bigg[x_j(t)-x_i(t)\bigg]\qquad i=1,2,\cdots,N \label{constantlin}
\end{eqnarray}
where $x_i(t)=(x_i^1(t),\cdots,x_i^n(t))^T\in R^n$;
$f(\cdot,t):R^n\times R^{+}\rightarrow R^n$ is continuous. The outer
coupling matrix $A$ satisfies: $a_{ij}\geq 0$, for $i\not=j$, and
let $a_{ii}=-\sum\limits_{j=1,j\not=i}^N a_{ij}$; while the inner
coupling matrix $\Gamma=diag(\gamma_1,\cdots,\gamma_n)$ is positive
definite. $c$ denotes the coupling strength.

If $\lim\limits_{t\rightarrow \infty}\|x_i(t)-x_j(t)\|=0$ for all
$i,j=1,2,\cdots,N$, where $\|\cdot\|$ denotes some norm, then
coupled systems (\ref{synchro}) are said to be synchronized
completely. In the following, we will investigate complete
synchronization. Hitherto, many approaches and criteria to ensure
synchronization have been derived (Lu et. al 2004 - Pikovsky,
Rosenblum, and Kurths 2001).

Moreover, in practice, the state variables $x_{i}(t)$ may be
unobservable; instead, we can observe $g(x_{i}(t))$, where
$g(\cdot)$ is a monotone increasing function. How to synchronize
coupled systems with the observed data is of great significance.
In (Chen, Zhu 2007), following nonlinearly coupled systems
\begin{eqnarray}
\dot{x}_i(t)=f(x_i(t),t)
+c\sum\limits_{j\not=i}a_{ij}\bigg[g(x_j(t))-g(x_i(t))\bigg] \qquad
i=1,2,\cdots,N\label{constantnonlin}
\end{eqnarray}
were proposed, where
$g(x_j(t))=(g_1(x_j^1(t)),\cdots,g_n(x_j^n(t))^T$ and every
$g_{i}(\cdot)$ is a nonlinear monotone increasing function.

It is pointed out in Pikovsky, Rosenblum, and Kurths 2001 that if
the coupling strength is larger than a critical value $c^{*}$,
then coupled systems (\ref{constantlin}) can be synchronized. In
fact, $c^{*}$ depends not only on the coupling matrix but also
depends on the dynamical behavior of function $f(x(t),t)$ in the
uncoupled system $\dot{x}(t)=f(x(t),t)$. Thus, the critical value
$c^{*}$ (suitable for all systems with various $f(x(t),t)$) is
much larger than the coupling strength $c$ needed for specified
coupled systems (see Lu and Chen 2004).

On the other hand, it is well known that a chaotic attractor
typically has embedded within it an infinite number of unstable
periodic orbits. In (Ott, Grebogi and Yorke 1990), it is shown
that one can convert a chaotic attractor to any one of a large
number of possible attracting time-periodic motions by making only
small time-dependent perturbations of an available system
parameter. On the other hand, if the attractor is not chaotic but
is, say, periodic, then small parameter perturbations can only
change the orbit slightly.

Therefore, how to synchronize a large number of specified chaotic
oscillators (for example, Lorenz oscillator or other oscillators)
with a relatively small coupling strength is an interesting problem.

For this purpose, we replace coupled systems (\ref{constantlin}) and
(\ref{constantnonlin}) by
\begin{eqnarray}
\dot{x}_i(t)=f(x_i(t),t)+c(t)\sum\limits_{j\not =i}a_{ij}\Gamma
\bigg[x_j(t)-x_i(t)\bigg]\qquad i=1,2,\cdots,N
\end{eqnarray}
\begin{eqnarray}
\dot{x}_i(t)=f(x_i(t),t)+c(t)\sum\limits_{j\not=i}
a_{ij}\bigg[g(x_j(t))-g(x_i(t))\bigg] \qquad i=1,2,\cdots,N
\end{eqnarray}
with an adaptive coupling strength $c(t)$.

The adaptation of parameters is widely used in the signal
processing and other research fields. Recently, the adaptive
approach is used to synchronize master-salve systems or mutually
coupled systems to a specified trajectory of the uncoupled node by
adding negative feedbacks (see Chen and Zhou 2006 - Boccaletti
2006). For example, in (Chen and Zhou 2006 and Zhou et.al 2006),
the authors investigated how to synchronize coupled systems:
\begin{eqnarray}
\dot{x}_i(t)=f(x_i(t),t)+h_i(x_1(t),x_2(t),\cdots,x_N(t)) \quad
1\leq i\leq N \label{MS}
\end{eqnarray}
to a specified trajectory $\dot{s}(t)=f(s(t),t) $ by adding controls
$-d_i(t)\bigg[x_i(t)-s(t)\bigg]$ with adaptation rule
$\dot{d}_i=k_i\|x_i(t)-s(t)\|_2^2$.

However, this adaptation rule does not suit complete
synchronization. Because, the synchronization trajectory, which
generally is not a trajectory of the uncoupled system, is unknown.

In this paper, new schemes to synchronize chaotic oscillators with
an adaptive coupling strength are proposed. It is revealed that a
large number of chaotic oscillators (including Lorenz oscillators,
Chen's oscillators, R\"ossler oscillators and Chua's circuits) can
be synchronized with a very small coupling strength. It means that
one can synchronize a large number of chaotic oscillators with very
weak coupling.

This paper is organized as follows: In section 2, some necessary
definitions, lemmas, and hypotheses are given; in section 3,
adaptive schemes for linearly and nonlinearly coupled system are
given; in section 4 and section 5, some simulations are given to
verify our theoretical results; while in section 6, an example is
given to illustrate that: for coupled periodic systems the coupling
strength should be very large; and the paper is concluded in section
7.

\section{Preliminaries}\quad
In this section, some definitions, denotations and lemmas throughout
the paper are presented.

\begin{mydef}\quad 
The set $\bf{S}$$=\{(x_1^T,x_2^T,\cdots,x_N^T)^T|x_i=x_j;\quad
i,j=1,2,\cdots,N\}$ is called the synchronization manifold.
\end{mydef}

\begin{mydef}\quad 
Matrix $A=(a_{ij})_{i,j=1}^N$  of order N is said to satisfy
Condition ${\bf A1}$, if

1. $a_{ij}\geq 0, i\not=j$, $a_{ii}=-\sum\limits_{j=1,j\not=i}^N
a_{ij},       i= 1,2, \cdots, N$

2. Eigenvalues of $A$ are all negative except an eigenvalue 0 with
multiplicity 1.

Furthermore, if $A\in A1$, and $a_{ij}=a_{ji}, i\not=j$, then we say
$A\in {\bf A2}$.
\end{mydef}

\begin{mydef}\quad 
Suppose that  $\Delta=diag\{\delta_1, \delta_2,\cdots,\delta_n\}$ is
a diagonal matrix, and $\varpi>0$. $f:R^n\times R^{+}\rightarrow
R^n$ is continuous. we say $f\in QUAD(\Delta,\varpi)$, if and only
if
\begin{equation}
(x-y)^T[f(x,t)-f(y,t)]-(x-y)^T\Delta(x-y)\leq
-\varpi(x-y)^T(x-y)\label{quadequ}
\end{equation}
holds for any $x,y\in R^n$.
\end{mydef}

\begin{myrem}\quad It is some control for the case
$(x-y)^T[f(x,t)-f(y,t)]>0$ (roughly speaking, $f(x,t)-f(y,t)$ and
$x-y$ have same sign), $||f(x,t)-f(y,t)||$ should be less than
$(x-y)^T[\Delta-\varpi I_{n}](x-y)$. Instead, if
$(x-y)^T[f(x,t)-f(y,t)]\le 0$, then inequality (\ref{quadequ}) is
satisfied automatically.

\end{myrem}

\begin{mylem}\quad 
If a matrix $A_{N\times N}\in{\bf A1}$. Then (see Wu 2005, Lu Chen
2006):

1. $[1, 1, \cdots, 1]^T$ is the right eigenvector of $A$.
corresponding to eigenvalue 0 with multiplicity 1;

2. The left eigenvector of $A$: $\xi=[\xi_{1}, \xi_2, \cdots,
\xi_N]^T \in R^N$ corresponding to eigenvalue 0 has following
properties: it is non-zero and its multiplicity is 1; all
$\xi_{i}\geq 0$, $i= 1, 2, \cdots, N.$ More precisely, $A$ is
irreducible if and only if all $\xi_{i}>0$, $i= 1, 2, \cdots, N$.
\end{mylem}

In the following, we always assume that {\bf $A$ is irreducible} and
$\sum\limits_{i=1}^N\xi_i=1$, $\xi_i>0$, which is called the
normalized left eigenvector. It is clear that if $A_{N\times N}\in
{\bf A2}$, then $\xi_{i}=\frac{1}{N}$ for $i=1,\cdots, N.$

\section{Main Results}\quad
In this section, we give the main results about how to design
synchronization algorithms with an adaptive coupling strength for
linearly and nonlinearly coupled systems.

\subsection{Linearly coupled systems}\quad
In this part, we discuss synchronization of linearly coupled systems
with the coupling strength adaptively.

\subsubsection{Constant coupling matrix}\quad
In this part, we investigate the case that the coupling matrix is
time independent.

Consider following linearly coupled systems with an adaptive
coupling strength:
\begin{equation}\label{lin}
\dot{x}_i(t)=f(x_i(t),t)+c(t)\sum\limits_{j\not=j} a_{ij}\Gamma
\bigg[x_j(t)-x_i(t)\bigg]
\end{equation}
the coupling matrix $A=(a_{ij})$ is not assumed to be symmetric.

For any irreducible (asymmetric) matrix $A_{N\times N}\in {\bf A1}$,
Denote $\Xi=diag\{\xi_1,\cdots,\xi_N\}$, $U=\Xi-\xi\xi^T$, where
$\xi=(\xi_1,\cdots,\xi_N)^T\in R^N$ is the normalized left
eigenvector corresponding to eigenvalue $0$.. It is easy to check
that $-U\in {\bf A2}$.

Furthermore, denote $X(t)=(x_1^T(t),\cdots,x_N^T(t))^T$,
$F(X(t))=(f(x_1(t),t)^T, \cdots, \\ f(x_N(t),t)^T)^T$,  ${\bf
\Delta}=I_N\otimes \Delta$, ${\bf A}=A\otimes \Gamma$, and ${\bf
U}=U \otimes I_n$, ${\bf \Xi}=\Xi \otimes I_n$, where $\otimes$ is
the Kronecker product. Then, (\ref{lin}) can be rewritten in a
compact form:
\begin{equation}
\dot{X}(t)=F(X(t))+c(t){\bf A}X(t)\label{system}
\end{equation}

The following Lemma gives a necessary and sufficient condition
whether a coupled system is synchronized. It is clear that the left
eigenvector $\xi$ plays a key role to prove synchronization.
\begin{mylem}\quad 
$X(t)\in{\bf S}$ if and only if
\begin{equation}
X^T(t){\bf U}X(t)=0
\end{equation}
\end{mylem}

In fact, for any $X(t)=(x_1^T(t),\cdots,x_N^T(t))^T$,
$Y(t)=(y_1^T(t),\cdots,y_N^T(t))^T$, we have (also see Wu 2005)
\begin{eqnarray}
X^T(t){\bf U}Y(t)=\frac{1}{2}\sum_{i,j=1}^{N}
\xi_{i}\xi_{j}(x_{i}(t)-x_{j}(t))^{T}(y_{i}(t)-y_{j}(t))\label{Chua}
\end{eqnarray}

Now, we propose a new scheme to synchronize linearly coupled systems
with an adaptive coupling strength and prove the following theorem.

\begin{mythm}\quad 
Suppose $A\in {\bf A1}$ is irreducible and $f\in
QUAD(\Delta,\varpi)$. Then, following coupled systems with an
adaptive coupling strength $c(t)$:
\begin{eqnarray}
\left\{\begin{array}{cc}\dot{X}(t)=F(X(t))+c(t){\bf A}X(t),\\
\dot{c}(t)=-\frac{\alpha}{2}X^T(t){\bf \Xi}{\bf A}X(t)
\end{array}\right. \label{Scheme}
\end{eqnarray}
can finally achieve synchronization with a relatively small coupling
strength, where $c(0)=0$ and $\alpha>0$.
\end{mythm}
{\bf Proof:}\quad Pick a sufficiently large constant $c>0$ and
define a Lyapunov function as:
\begin{equation}
V(X(t))=\frac{1}{2}X^T(t){\bf U}X(t)+\frac{1}{\alpha}(c-c(t))^2
\end{equation}
Then, noting $UA=\Xi A$, we have
\begin{eqnarray*}
\frac{dV(X(t))}{dt}&=&X(t)^T{\bf U}[F(X(t))+c(t){\bf
A}X(t)]+(c-c(t))X^T(t){\bf \Xi}{\bf A}X(t)\\
&=&X(t)^T{\bf U}[F(X(t))-{\bf \Delta} X(t)]+X^T(t){\bf U}{\bf
\Delta}X(t)+cX^T(t){\bf \Xi}{\bf A}X(t)
\end{eqnarray*}
By the identity (\ref{Chua}) and the assumption $f\in
QUAD(\Delta,\varpi)$, we have
\begin{eqnarray*}
&&X^T(t){\bf U}[F(X(t))-{\bf \Delta} X(t)]\\
&=& \sum_{i,j=1;i\not=j}^{N}
\xi_{i}\xi_{j}[x_{i}(t)-x_{j}(t)]^{T}[f(x_{i}(t),t)-f(x_{j}(t),t)
-\Delta(x_{i}(t)-x_{j}(t))]\\
&\le& -\varpi\sum_{i,j=1;i\not=j}^{N}
\xi_{i}\xi_{j}[x_{i}(t)-x_{j}(t)]^{T}[x_{i}(t)-x_{j}(t)]=-\varpi
X^T(t){\bf U}X(t)
\end{eqnarray*}
Therefore,
\begin{eqnarray}
\frac{dV(t)}{dt} \leq -\varpi X(t)^T{\bf U}X(t)+X^T(t){\bf U}{\bf
\Delta}X(t)+cX^T(t){\bf \Xi}{\bf A}X(t)\label{a}
\end{eqnarray}
Write $\tilde{x}_j(t)=(x_1^j(t),\cdots,x_N^j(t))^T$ for
$j=1,2,\cdots,n$, we have
\begin{eqnarray}
X^T(t){\bf U}({\bf \Delta}+c{\bf A})X(t)=\sum\limits_{j=1}^n
\tilde{x}_j^T(t)\delta_j U\tilde{x}_j(t)+c\sum\limits_{j=1}^n
\gamma_j \tilde{x}_j^T(t)\Xi A\tilde{x}_j(t)\label{a1}
\end{eqnarray}

It can be seen that $\Xi A+A^{T}\Xi$ is a symmetric matrix with
negative diagonal and row-sum zero. Now, let $v_{1},\cdots,v_{N}$ be
the normalized eigenvectors of the matrix $\frac{1}{2}(\Xi
A+A^{T}\Xi)$ with corresponding eigenvalues
$\lambda_{1}=0>\lambda_{2}\geq\cdots\geq\lambda_{N}$. Moreover, by
$\tilde{x}^{*}_j(t)$ denote the projection of $\tilde{x}_j(t)$ on
the subspace $L$ spanned by $v_{2},\cdots,v_{N}$. Then,
\begin{eqnarray}
&&c\sum\limits_{j=1}^n \tilde{x}_j^T(t)\gamma_j \Xi
A\tilde{x}_j(t)=\frac{c}{2}\sum\limits_{j=1}^n
\tilde{x}_{j}^{{T}}(t)\gamma_j (\Xi
A+A^{T}\Xi)\tilde{x}_j(t)\nonumber\\&&=\frac{c}{2}\sum\limits_{j=1}^n
\tilde{x}_{j}^{{*}^{T}}(t)\gamma_j (\Xi
A+A^{T}\Xi)\tilde{x}^{*}_j(t)\le c\lambda_{2}\sum\limits_{j=1}^n
\gamma_j\tilde{x}_{j}^{{*}^{T}}(t) \tilde{x}^{*}_j(t)
\end{eqnarray}
\begin{eqnarray}
\sum\limits_{j=1}^n \tilde{x}_j^T(t)\delta_j
U\tilde{x}_j(t)=\sum\limits_{j=1}^n
\tilde{x}_{j}^{{*}^{T}}(t)\delta_j U\tilde{x}^{*}_j(t)
\end{eqnarray}
Substituting into (\ref{a}), we have
\begin{eqnarray}
\frac{dV(t)}{dt}\le -\varpi X(t)^T{\bf U}X(t)+\sum\limits_{j=1}^n
\tilde{x}_{j}^{{*}^{T}}(t)\delta_j
U\tilde{x}^{*}_j(t)+c\lambda_{2}\sum\limits_{j=1}^n
\gamma_j\tilde{x}_{j}^{{*}^{T}}(t) \tilde{x}^{*}_j(t)
\end{eqnarray}
It is clear that if $c$ is sufficient large, then
\begin{eqnarray}
\sum\limits_{j=1}^n \tilde{x}_{j}^{{*}^{T}}(t)\delta_j
U\tilde{x}^{*}_j(t)+c\lambda_{2}\sum\limits_{j=1}^n
\gamma_j\tilde{x}_{j}^{{*}^{T}}(t) \tilde{x}^{*}_j(t)<0
\end{eqnarray}
In summary, we have $\frac{dV(t)}{dt}\leq-\varpi X(t)^T{\bf
U}X(t)\leq 0 $.

Similarly, it can be seen that $\dot{c}(t)\geq 0$ and $\dot{c}(t)=0$
if and only if $X\in {\bf S}$, so $c(t)>0$ for $t>0$.

It is obvious that $\dot{V}=0$ if and only if $X\in {\bf S}$.
According to the well-known Lyapunov-LaSall type theorem for
functional differential equations (see Kuang 1993), the trajectory
of coupled systems, starting with arbitrary initial value,
converges asymptotically to the largest invariant set $H_1$
contained in $H_2=\{\dot{V}(t)=0\}$ as $t\rightarrow +\infty$,
where $H_1=\{[X,c]^T:X^{T}UX=0, c=c_0\in R^{+}\}$. Therefore,
$X^{T}(t)UX(t)\rightarrow 0$ and $c(t)\rightarrow c_0$ for some
constant $c_0$.

Theorem 1 is proved completely.

\begin{myrem}\quad
Adaptive algorithm is often used in many research fields. Many
authors use $-d_i\bigg[x_i(t)-s(t)\bigg]$ as the negative feedback
adaptation, where $\dot{d}_i=k_i\|x_i(t)-s(t)\|_2^2$, $k_i$ are
positive constants and $s(t)$ is a special solution (or an
equilibrium point) of the uncoupled system. Here, the
synchronization state is unknown. Therefore, previous negative
feedback adaptation methods are invalid. Theorem 1 provides a new
adaptive algorithm, which succeeds in the synchronization of complex
networks with an adaptive coupling strength.
\end{myrem}

\subsubsection{Unknown constant coupling matrix}\quad
In some physical coupled systems, the coupling matrix  may be
unknown, though we know $A\in A1$. Can we design an adaptive
algorithm to synchronize coupled systems with an adaptive coupling
strength for a unknown coupling matrix? In this subsection, we will
give an affirmative answer by proving the following theorem.

\begin{mythm}\quad 
Suppose that the unknown coupling matrix $A\in {\bf A1}$ is
irreducible and $f\in QUAD(\Delta,\varpi)$. Then, following coupled
systems with an adaptive coupling strength $c(t)$:
\begin{eqnarray}
\left\{\begin{array}{cc}\dot{X}(t)=F(X(t))+c(t){\bf A}X(t),\\
\dot{c}(t)=-\frac{\alpha}{2}X^T(t){\bf \tilde{A}}X(t)
\end{array}\right.
\end{eqnarray}
can achieve synchronization with a relatively small coupling
strength, where ${\bf \tilde{A}}=\tilde{A}\otimes I_n$,
$\tilde{A}\in {\bf A2}$ is any irreducible matrix, $c(0)=0$ and
$\alpha>0$.
\end{mythm}
{\bf Proof:}\quad Define a slightly different Lyapunov function
\begin{equation}
V_{1}(X(t))=\frac{1}{2}X^T(t){\bf
U}X(t)+\frac{\eta}{\alpha}(c-c(t))^2
\end{equation}
where $\eta>0$ is a constant. Then, noting $UA=\Xi A$, we have
\begin{eqnarray*}
\frac{dV_{1}(X(t))}{dt}&=&X(t)^T{\bf U}[F(X(t))+c(t){\bf
A}X(t)]+\eta (c-c(t))X^T(t){\bf \tilde{A}}X(t)\\
&=&X(t)^T{\bf U}[F(X(t))-{\bf \Delta} X(t)]+X^T(t){\bf U}{\bf
\Delta}X(t)\\
&&+c(t)X^T(t){\bf \Xi}{\bf A}X(t)+\eta (c-c(t))X^T(t){\bf \tilde{A}}X(t)\\
&\leq&-\varpi X(t)^T{\bf U}X(t)+X^T(t)\bigg({\bf U}{\bf
\Delta}+c\eta{\bf \tilde{A}}\bigg)X(t)\\
&&+c(t)X^T(t)\bigg({\bf \Xi}{\bf A}-\eta{\bf \tilde{A}}\bigg)X(t)
\end{eqnarray*}

Now, pick a sufficiently large constant $c$ and a sufficiently small
$\eta$ such that
\begin{eqnarray}
\left\{
\begin{array}{c}
{\bf U}{\bf \Delta}+c\eta{\bf \tilde{A}}\le 0\\
\{{\bf \Xi}{\bf A}\}^s-\eta{\bf \tilde{A}}\le 0
\end{array}
\right.
\end{eqnarray}
Therefore, $\frac{dV_{1}(X(t))}{dt}\leq -\varpi X(t)^T{\bf U}X(t)\le
0$. By similar arguments used in Theorem 1, we can prove Theorem 2.

\subsubsection{Time-varying coupling
matrix}\quad In practice, the coupling matrix is time-dependent,
which means that the coupling matrix changes along with time.
Therefore, it is natural to investigate following linearly coupled
systems with a time-varying coupling matrix:
\begin{equation}
\dot{X}(t)=F(X(t))+c(t){\bf A}(t)X(t)\label{system1}
\end{equation}

The following theorem for time-varying coupling can be proved with
the same arguments used in the proof of Theorem 1.
\begin{mythm}\quad 
Suppose that $f\in QUAD(\Delta,\varpi)$, the coupling matrix
$A(t)\in {\bf A1}$ is irreducible and has the same left eigenvector
$\xi$ corresponding to eigenvalue $0$ for all $t$. If the largest
non-zero eigenvalue of the matrix $\Xi A(t)+A^{T}(t)\Xi$ satisfies
$\lambda_{2}(t)\le \lambda<0$ for all $t$. Then, following coupled
systems with an adaptive coupling strength $c(t)$:
\begin{eqnarray}
\left\{\begin{array}{cc}\dot{X}(t)=F(X(t))+c(t){\bf A(t)}X(t),\\
\dot{c}(t)=-\frac{\alpha}{2}X^T(t){\bf \Xi}{\bf A(t)}X(t)
\end{array}\right.\label{time-varying}
\end{eqnarray}
can finally achieve synchronization with a relatively small coupling
strength, where $c(0)=0$ and $\alpha>0$.
\end{mythm}

\begin{myrem}\quad  
The condition that $A(t)\in {\bf A1}$ is irreducible and has the
same left eigenvector $\xi$ corresponding to eigenvalue $0$ for
all $t$ looks quite strong. However, it just suits the
node-balanced coupling networks discussed in (Belykh, V. et.al
2006), where the coupling matrix $A(t)$ is assumed to be row-sum
zero as well as column-sum zero for each $t$. In this case,
$A(t)\in {\bf A1}$ has the same left eigenvector $[1,\cdots,1]^T$
corresponding to eigenvalue $0$ for all $t$. Therefore, Theorem 3
applies to the node-balanced coupling networks.
\end{myrem}

As a direct consequence of Theorem 3, we can obtain the following
simple adaptive scheme, if every $A(t)$ is symmetric.

\begin{mycol}\quad 
Suppose that $f\in QUAD(\Delta,\varpi)$, the coupling matrix
$A(t)\in {\bf A2}$ is irreducible. If the largest non-zero
eigenvalue of the matrix $A(t)$ satisfies $\lambda_{2}(t)\le
\lambda<0$, for all $t$. Then, following coupled systems with an
adaptive coupling strength $c(t)$:
\begin{eqnarray}
\left\{\begin{array}{cc}\dot{X}(t)=F(X(t))+c(t){\bf A(t)}X(t),\\
\dot{c}(t)=-\frac{\alpha}{2}X^T(t){\bf A(t)}X(t)
\end{array}\right.
\end{eqnarray}
can finally achieve synchronization with a relatively small coupling
strength, where $c(0)=0$ and $\alpha>0$.
\end{mycol}

\begin{myrem}\quad 
In Theorem 3 and Corollary 1, condition $\lambda_{2}(t)\le
\lambda<0$ plays a key role. However, calculating $\lambda_{2}(t)$
for all $t$ is impossible numerically. If all $A(t)$ can be
dominated by a constant matrix $\hat{A}$, then we can synchronize
chaotic oscillators directly using (\ref{Scheme}), only replacing
$A$ in the adaptive algorithm by $\hat{A}$. Following corollary
explains it.
\end{myrem}

\begin{mycol}\quad 
Suppose that $f\in QUAD(\Delta,\varpi)$, the coupling matrix
$A(t)\in {\bf A2}$ is irreducible. If there exists a constant matrix
$\hat{A}=(\hat{a}_{ij})\in A2$ such that $\hat{a}_{ij}\geq
{a}_{ij}(t), i\not=j$ hold for all $t$. Then, following coupled
systems with an adaptive coupling strength $c(t)$
\begin{eqnarray}
\left\{\begin{array}{cc}\dot{X}(t)=F(X(t))+c(t){\bf A(t)}X(t),\\
\dot{c}(t)=-\frac{\alpha}{2}X^T(t){\bf \hat{A}}X(t)
\end{array}\right.
\end{eqnarray}
can finally achieve synchronization with a relatively small coupling
strength, where $c(0)=0$ and $\alpha>0$.
\end{mycol}
{\bf Proof}\quad In this case, $\xi_{i}=\frac{1}{N}$, which implies
$\Xi A(t)=\frac{1}{N}A(t)$, and
\begin{eqnarray}
&&X^T(t){\bf A}(t)X(t)=\sum_{i,j=1;i\not=j}^{N}
a_{ij}(t)(x_{i}(t)-x_{j}(t))^{T}(x_{i}(t)-x_{j}(t))
\nonumber\\&&\le\sum_{i,j=1;i\not=j}^{N}
\hat{a}_{ij}(x_{i}(t)-x_{j}(t))^{T}(x_{i}(t)-x_{j}(t))=X^T(t){\bf
\hat{A}}X(t)
\end{eqnarray}
Using the Lyapunov function
\begin{equation}
V(X(t))=\frac{1}{2}X^T(t){\bf U}X(t)+\frac{1}{\alpha}(c-c(t))^2
\end{equation}
we have
\begin{eqnarray*}
\frac{dV(t)}{dt} &\le &X(t)^T{\bf U}[F(X(t))-{\bf \Delta}
X(t)]+X^T(t){\bf U}({\bf \Delta}+cN{\bf \hat{A}})X(t)
\end{eqnarray*}
By the same arguments used in the proof of Theorem 1. we can obtain
Corollary 2.

\subsection{Nonlinearly coupled systems}\quad
We consider following nonlinearly coupled systems with an adaptive
coupling strength:
\begin{equation}\label{nonlin}
\dot{x}_i(t)=f(x_i(t),t)
+c(t)\sum\limits_{j\not=i}a_{ij}\bigg[g(x_j(t))-g(x_i(t))\bigg]
\end{equation}
where $g(x_i(t))=[g_{1}(x_i^1(t)),\cdots,g_{n}(x_i^n(t))]^{T}$.

Denote $G(X(t))=(g^T(x_1(t)),\cdots,g^T(x_N(t)))^T$ and ${\bf
A}=A\otimes I_n$, then we have
\begin{equation}
\frac{dX(t)}{dt}=F(X(t))+c(t){\bf A}G(X(t))
\end{equation}

\begin{mythm}\quad 
Suppose $A\in {\bf A2}$ is irreducible, $f\in QUAD(\Delta,\varpi)$,
$\frac{g_k(u)-g_k(v)}{u-v}\geq \beta$ for $\forall u\not= v,
k=1,\cdots,n$, where $\beta$ is a positive constant. Then, following
nonlinearly coupled systems with an adaptive coupling strength
$c(t)$:
\begin{eqnarray}
\left\{\begin{array}{cc}\dot{X}(t)=F(X(t))+c(t){\bf A}G(X(t)),\\
\dot{c}(t)=-\frac{\alpha}{2}X^T(t){\bf A}G(X(t))
\end{array}\right.
\end{eqnarray}
will be synchronized with a relatively small coupling strength,
where $c(0)=0$ and $\alpha>0$.
\end{mythm}

{\bf Proof:}\quad Because $A\in {\bf A2}$, then
$\xi_{i}=\frac{1}{N}$ and $UA=\frac{1}{N}A$. Using the same Lyapunov
function as in the Theorem 1, we have
\begin{eqnarray*}
\frac{dV(t)}{dt}&=&X(t)^T{\bf U}[F(X(t))+c(t){\bf
A}G(X(t))]+(c-c(t))X^T(t){\bf U}{\bf A}G(X(t))\\
&=&X(t)^T{\bf U}[F(X(t))-{\bf \Delta} X(t)]+X^T(t){\bf U}{\bf
\Delta}X(t)+cNX^T(t){\bf U}{\bf A}G(X(t))
\end{eqnarray*}
Furthermore,
\begin{eqnarray*}
&&X^T(t){\bf U}{\bf A}G(X(t))=\sum\limits_{i=1}^n
\tilde{x}_i^T(t)UA\tilde{g}_{i}(\tilde{x}_i(t))
=\sum\limits_{i=1}^n\frac{1}{N}\tilde{x}_i(t)^TA\tilde{g}_{i}(\tilde{x}_i(t))\\
&=&-\sum\limits_{i=1}^n\frac{1}{N}\sum\limits_{j>k}a_{jk}
[x^j_i(t)-x^k_i(t)]^T[g_i(x^j_i(t))
-g_i(x^k_i(t))]\\
&\leq&-\beta\sum\limits_{i=1}^n\frac{1}{N}\sum\limits_{j>k}a_{jk}[x^j_i(t)-x^k_i(t)]^T
[x^j_i(t)-x^k_i(t)]\\
&=&\beta\sum\limits_{i=1}^n\frac{1}{N}\tilde{x}_i(t)^TA\tilde{x}_i(t)
=\beta\sum\limits_{i=1}^n\tilde{x}_i(t)^TUA\tilde{x}_i(t)
\end{eqnarray*}
where
$\tilde{g}_{i}(\tilde{x}_i(t))=(g_i(x_i^1(t)),\cdots,g_i(x_i^N(t)))^T$,
for $i=1,2,\cdots,n$. Therefore,
\begin{equation}
\frac{dV(t)}{dt}\leq -\varpi X(t)^T{\bf U}X(t)+\sum\limits_{j=1}^n
\tilde{x}_j^T(t)[\delta_j U+c\beta N UA]\tilde{x}_j(t)
\end{equation}
By similar arguments used in Theorem 1, we complete the proof of
Theorem 4.

Similar to linearly coupled systems, for nonlinearly coupled systems
with a time-varying coupling matrix:
\begin{eqnarray}
\frac{dX(t)}{dt}=F(X(t))+c(t){\bf A}(t)G(X(t))
\end{eqnarray}
we have
\begin{mycol}\quad 
Suppose $A(t)\in {\bf A2}$ is irreducible, $f\in
QUAD(\Delta,\varpi)$, $\frac{g_k(u)-g_k(v)}{u-v}\geq \beta$ for
$\forall u\not= v, k=1,\cdots,n$, where $\beta$ is a positive
constant. If the largest non-zero eigenvalue of $A(t)$ satisfies
$\lambda_{2}(t)\le \lambda<0$, for all $t$. Then, following
nonlinearly coupled systems with an adaptive coupling strength
\begin{eqnarray}
\left\{\begin{array}{cc}\dot{X}(t)=F(X(t))+c(t){\bf A}(t)G(X(t)),\\
\dot{c}(t)=-\frac{\alpha}{2}X^T(t){\bf A}(t)G(X(t))
\end{array}\right.
\end{eqnarray}
will be synchronized with a relatively small coupling strength,
where $c(0)=0$ and $\alpha>0$.
\end{mycol}

\section{Numerical Simulation}\quad
To validate the effectiveness of the proposed synchronization
algorithms with an adaptive coupling strength, in the following, we
couple 100 chaotic oscillators by following linearly coupled systems
with adaptive strength
\begin{eqnarray}
\left\{\begin{array}{cc}\dot{X}(t)=F(X(t))+c(t){\bf A}X(t),\\
\dot{c}(t)=-\frac{\alpha}{2}X^T(t){\bf \Xi}{\bf A}X(t)
\end{array}\right. \label{equ1}
\end{eqnarray}
and nonlinearly coupled systems with adaptive strength
\begin{eqnarray}
\left\{\begin{array}{cc}\dot{X}(t)=F(X(t))+c(t){\bf A}G(X(t)),\\
\dot{c}(t)=-\frac{\alpha}{2}X^T(t){\bf B}G(X(t))
\end{array}\right. \label{equ2}
\end{eqnarray}
where $X(t)=(x_1^T(t),\cdots,x_{100}^T(t))^T$,
$F(X(t))=(f(x_1(t),t)^T,\cdots,f(x_{100}(t),t)^T)^T$, ${\bf
A}=A\otimes I_{3}$, $A\in R^{100\times 100}$. ${\bf B}=B\otimes
I_{3}$, and $B=\frac{1}{2}(A+A^{T})$,
$G(X(t))=(g(x_1(t)),g(x_2(t)),\cdots,\newline g(x_{100}(t)))^T$,
$g(x_i(t))=(x_i^1(t)+tanh(x_i^1(t)),x_i^2(t)+tanh(x_i^2(t)),
x_i^3(t)+tanh(x_i^3(t)))^T$, $i=1,\cdots,100$. The prototypes of
$f(\cdot)$ are Chua's circuit, Chen's oscillator, Lorenz's
oscillator and R\"ossler's oscillator, respectively. For the
coupling matrix $A$, first, we construct a coupling matrix of a
small-world network generated by the method proposed in (Watts
et.al 1998). Then, replace every non-zero element $a_{ij},
i\not=j$ with a positive random scalar, which is equally
distributed in $[0,1]$, and choose the diagonal elements to ensure
$A\in {\bf A1}$.

We also couple 100 chaotic oscillators by following linearly coupled
systems with an adaptive strength
\begin{eqnarray}
\left\{\begin{array}{cc}\dot{X}(t)=F(X(t))+c(t){\bf A}X(t),\\
\dot{c}(t)=-\frac{\alpha}{2}X^T(t){\bf \tilde{A}}X(t)
\end{array}\right. \label{any}
\end{eqnarray}
where ${\bf \tilde{A}}=\tilde{A}\otimes I_{3}$, $\tilde{A}\in
R^{100\times 100}$ are globally connected matrix or randomly
connected matrix respectively.

In all simulations, we use the quantity
$E(t)=\sqrt{\sum\limits_{j=2}^{100}\|x_i(t)-x_1(t)\|^2/99}$ as a
measure of synchronization error.

\subsection{Chua's Circuits}

The uncoupled equation is :
\begin{eqnarray*}
\left \{
\begin{array}{l}
\frac{dx}{dt}=\bar{m}[y-\bar{h}(x)]\\
\frac{dy}{dt}=x-y+z\\
\frac{dz}{dt}=-\bar{n}y
\end{array}
\right .
\end{eqnarray*}
where $\bar{h}(x)=\frac{2}{7}x-\frac{3}{14}[|x+1|-|x-1|]$,
$\bar{m}=9$ and $\bar{n}=14\frac{2}{7}$.  Figure 1. (a) and (b) show
the dynamics of $c(t)$ and $E(t)$ for linearly coupled systems
(\ref{equ1}); while (c) and (d) show the dynamics of $c(t)$ and
$E(t)$ for nonlinearly coupled systems (\ref{equ2}).
\begin{center}
\includegraphics[height=4.5in,width=5.5in]{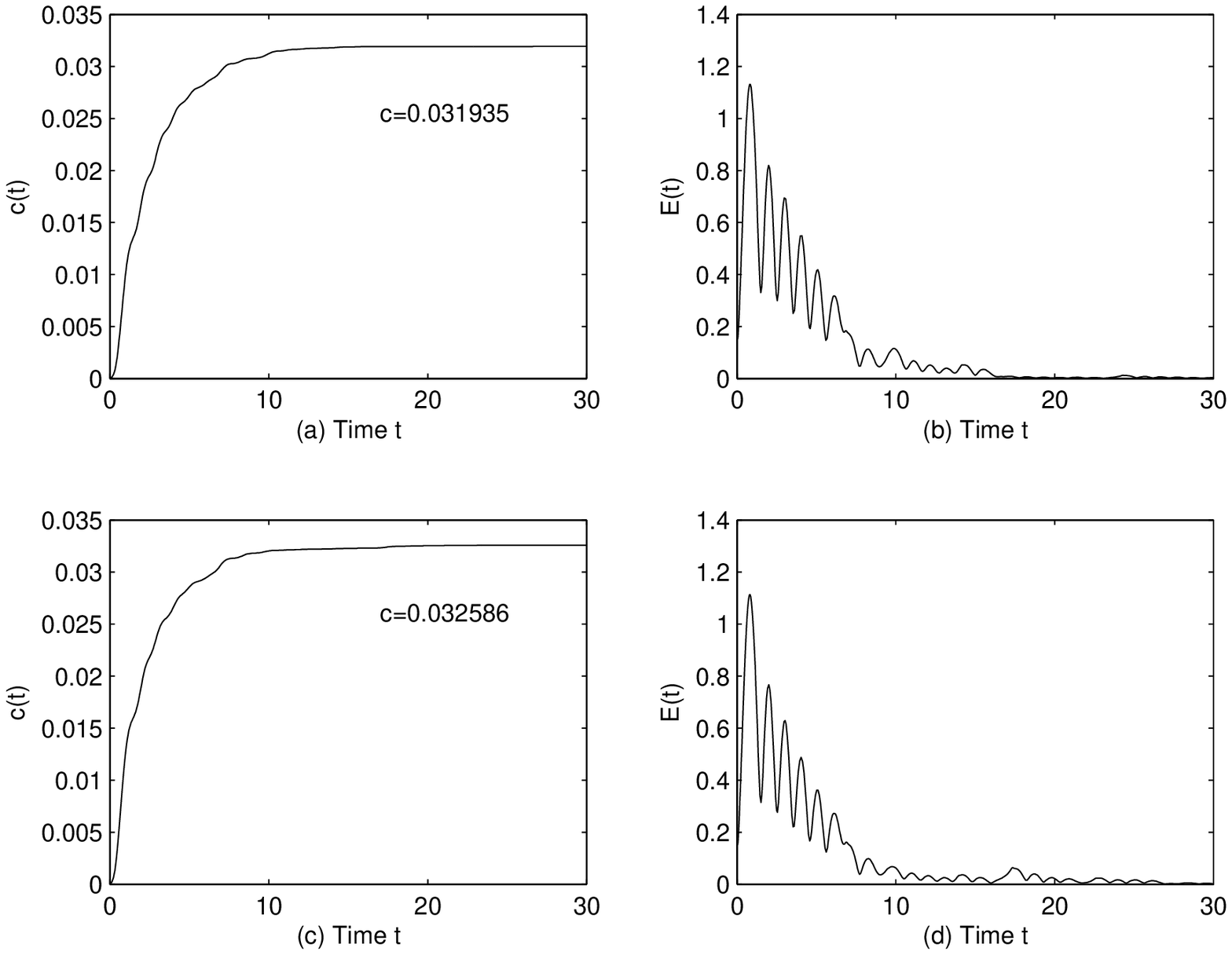}
\end{center}
\begin{center}
Figure 1. The dynamics of $c(t)$ and $E(t)$ for 100 linearly and
nonlinearly\\
coupled chua's circuits with an adaptive coupling strength
\end{center}

Figure 2. (a) and (b) show the dynamics of $c(t)$ and $E(t)$ for
linearly coupling systems (\ref{any}) when $\tilde{A}$ is the
globally coupled matrix; while (c) and (d) show the dynamics of
$c(t)$ and $E(t)$ for linearly coupled systems (\ref{any}) when
$\tilde{A}\in {\bf A2}$ is a random matrix.
\begin{center}
\includegraphics[height=4.5in,width=5.5in]{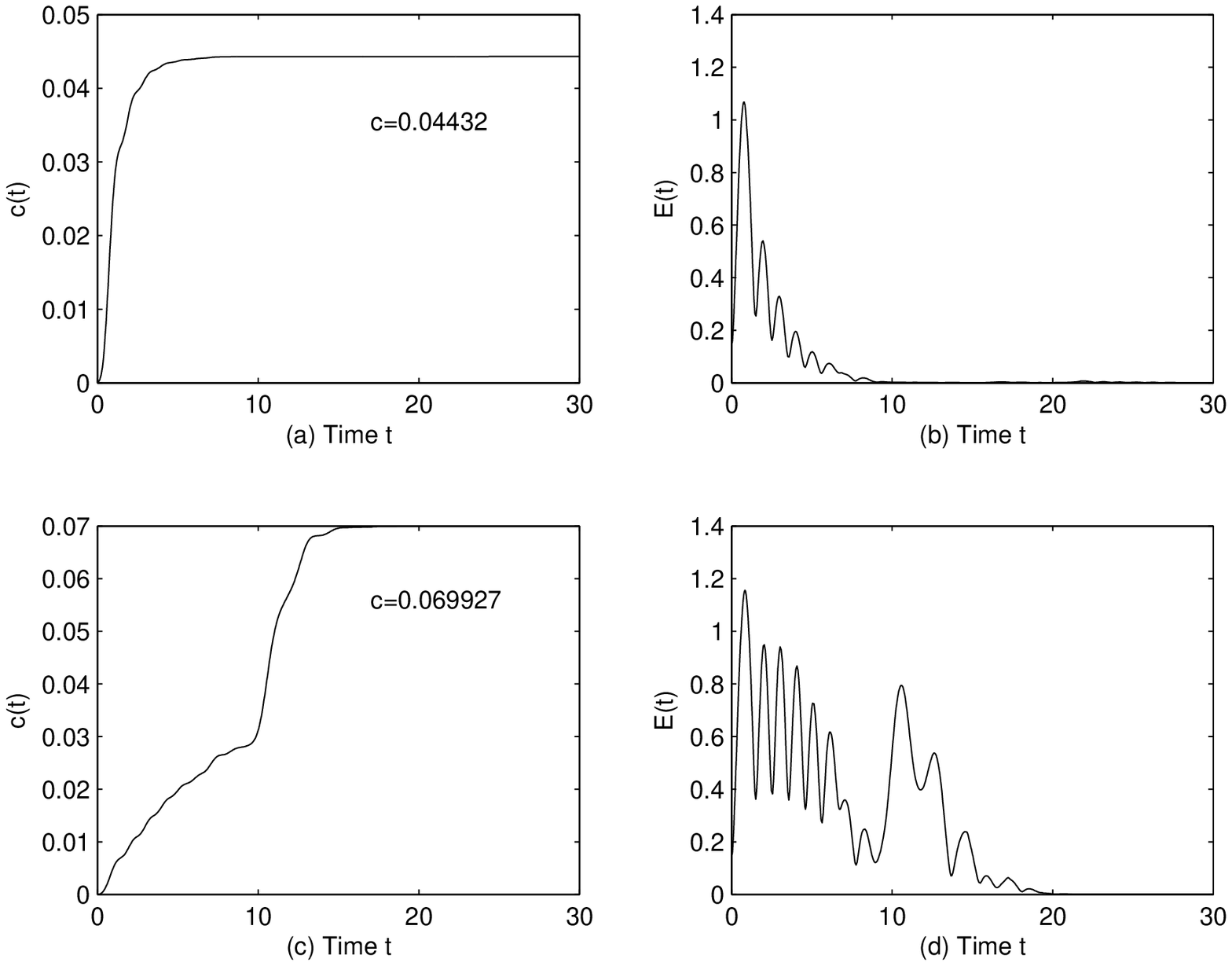}
\end{center}
\begin{center}
Figure 2. The dynamics of $c(t)$ and $E(t)$ for (\ref{any}) with
globally and randomly connected $\tilde{A}$
\end{center}
\subsection{Chen's Oscillator}

The uncoupled equation is :
\begin{eqnarray*}
\left\{
\begin{array}{l}
\dot{x}=a(y-x)\\
\dot{y}=(c-a)x-xz+cy\\
\dot{z}=xy-bz
\end{array}
\right.
\end{eqnarray*}
where $a=35$, $b=3$ and $c=28$. Figure 3. (a) and (b) show the
dynamics of $c(t)$ and $E(t)$ for linearly coupled systems
(\ref{equ1}); while (c) and (d) show the dynamics of $c(t)$ and
$E(t)$ for nonlinearly coupled systems (\ref{equ2}).

\begin{center}
\includegraphics[height=4.5in, width=5.5in]{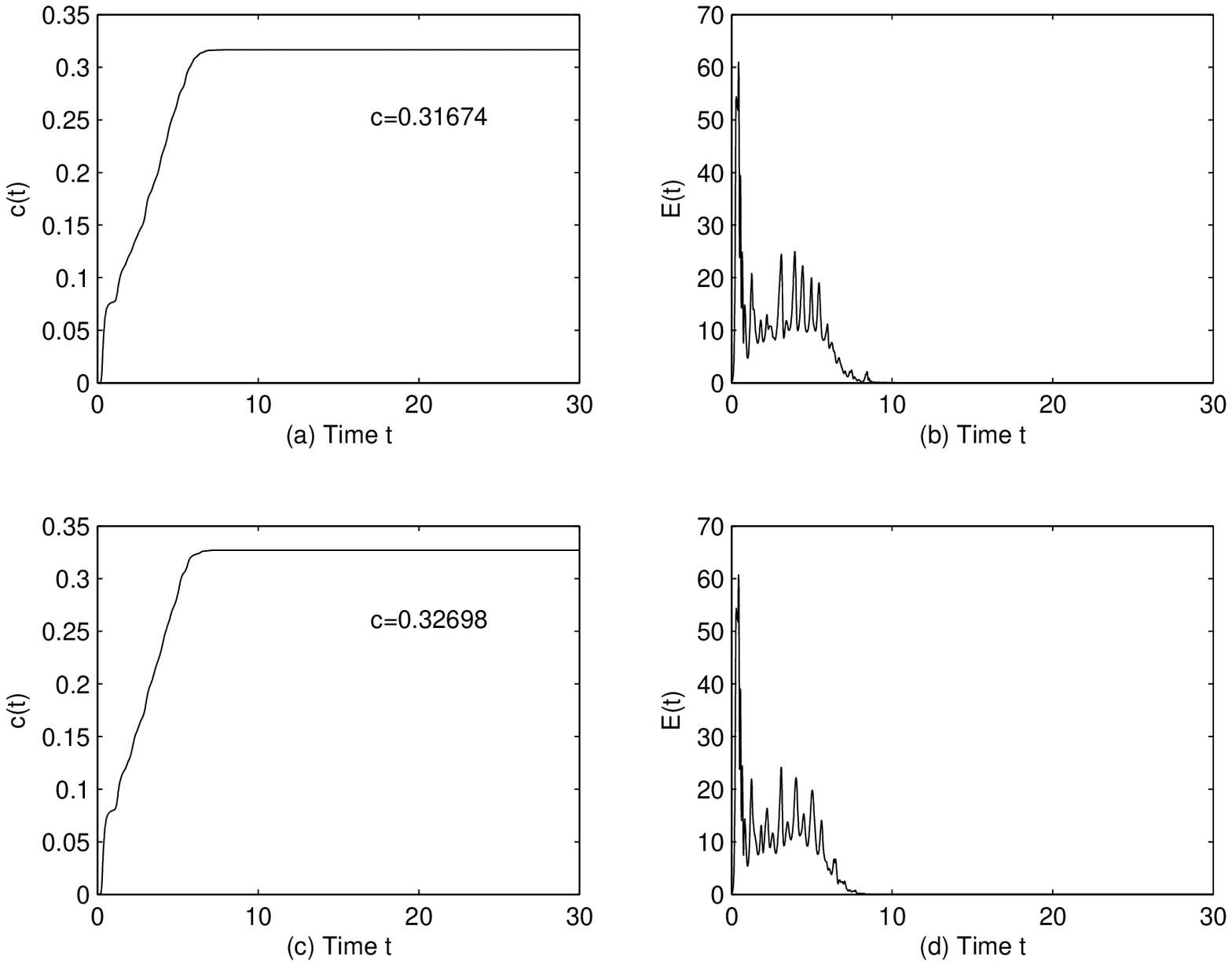}
\end{center}
\begin{center}
Figure 3. The dynamics of $c(t)$ and $E(t)$ for 100 linearly and
nonlinearly\\
coupled Chen's chaotic oscillators with an adaptive coupling
strength
\end{center}

Figure 4. (a) and (b) show the dynamics of $c(t)$ and $E(t)$ for
linearly coupling systems (\ref{any}) when $\tilde{A}$ is the
globally coupled matrix; while (c) and (d) show the dynamics of
$c(t)$ and $E(t)$ for linearly coupled systems (\ref{any}) when
$\tilde{A}\in {\bf A2}$ is a random matrix.
\begin{center}
\includegraphics[height=4.5in,width=5.5in]{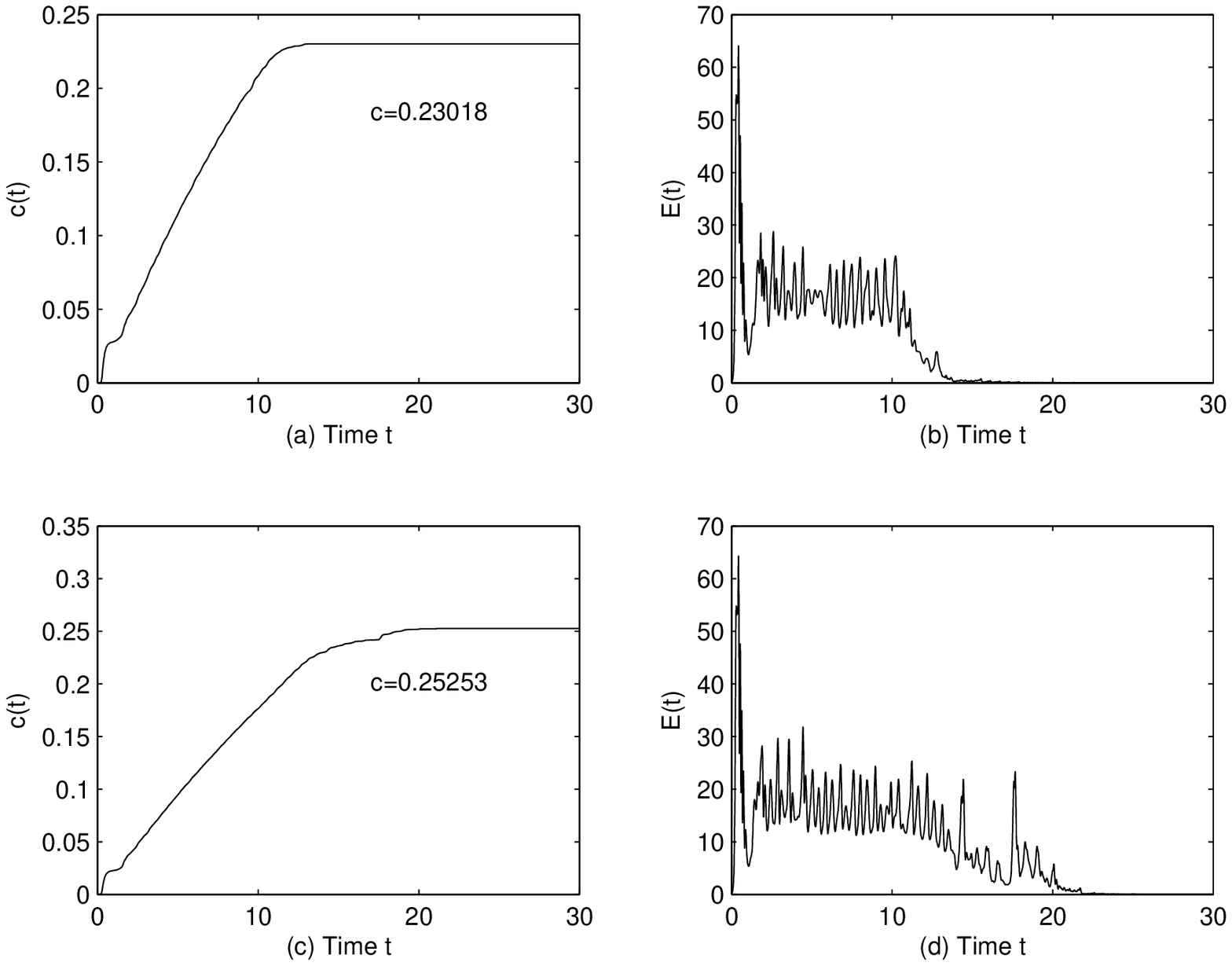}
\end{center}
\begin{center}
Figure 4. The dynamics of $c(t)$ and $E(t)$ for (\ref{any}) with
globally and randomly connected $\tilde{A}$
\end{center}
\subsection{Lorenz's Oscillator}

The uncoupled equation is :
\begin{eqnarray*}
\left\{
\begin{array}{l}
\dot{x}_1=\beta (x_2-x_1)\\
\dot{x}_2=\alpha x_1 -x_1 x_3 -x_2\\
\dot{x}_3=x_1 x_2-b x_3
\end{array}
\right.
\end{eqnarray*}
where $\beta=10$, $\alpha=28$, and $b=\frac{8}{3}$. Figure 5. (a)
and (b) show the dynamics of $c(t)$ and $E(t)$ for linearly coupled
systems (\ref{equ1}); while (c) and (d) show the dynamics of $c(t)$
and $E(t)$ for nonlinearly coupled systems (\ref{equ2}).

\begin{center}
\includegraphics[height=4.5in, width=5.5in]{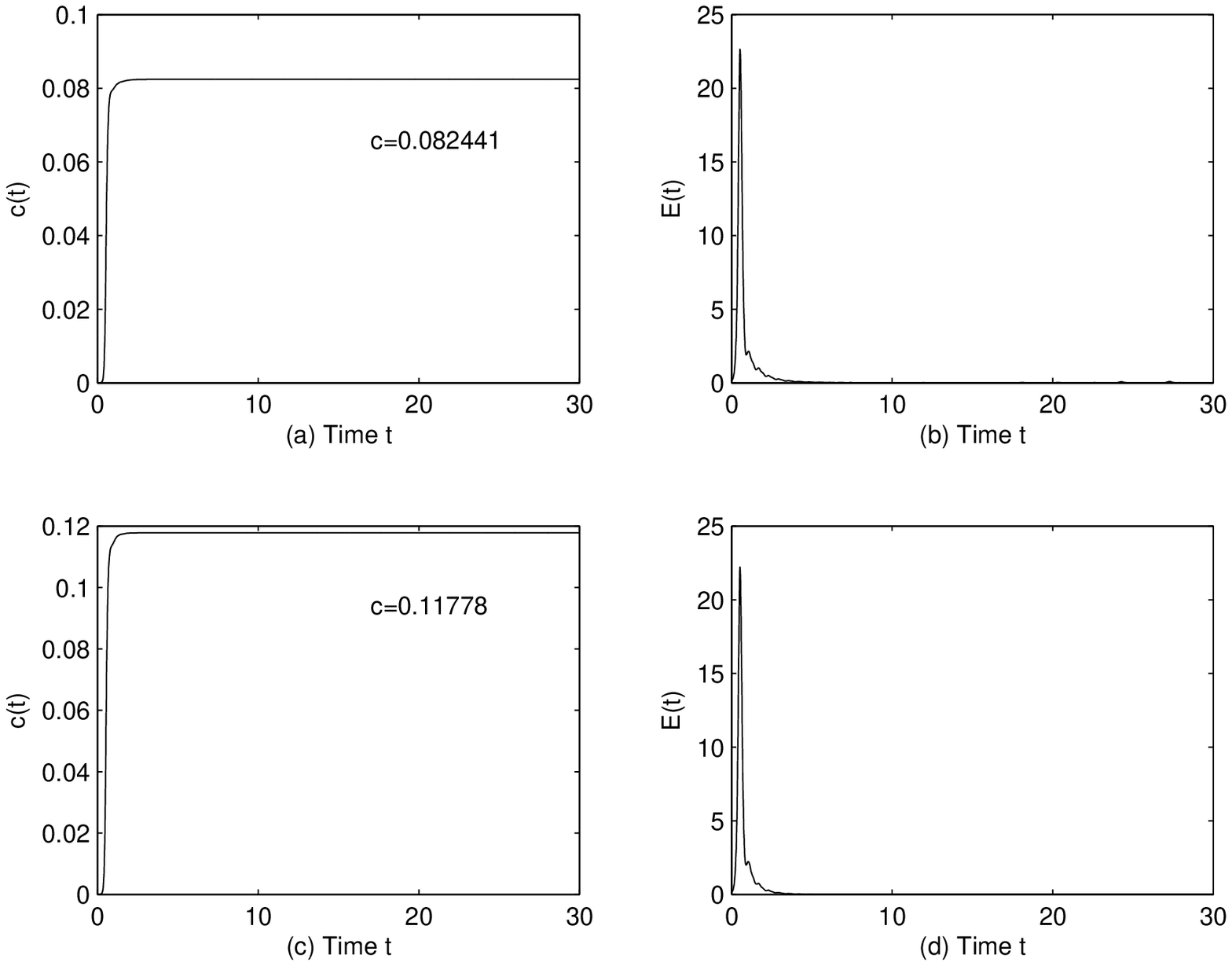}
\end{center}
\begin{center}
Figure 5. The dynamics of $c(t)$ and $E(t)$ for 100 linearly and
nonlinearly \\
coupled Lorenz's chaotic oscillators with an adaptive coupling
strength
\end{center}

Figure 6. (a) and (b) show the dynamics of $c(t)$ and $E(t)$ for
linearly coupling systems (\ref{any}) when $\tilde{A}$ is the
globally coupled matrix; while (c) and (d) show the dynamics of
$c(t)$ and $E(t)$ for linearly coupled systems (\ref{any}) when
$\tilde{A}\in {\bf A2}$ is a random matrix.
\begin{center}
\includegraphics[height=4.5in,width=5.5in]{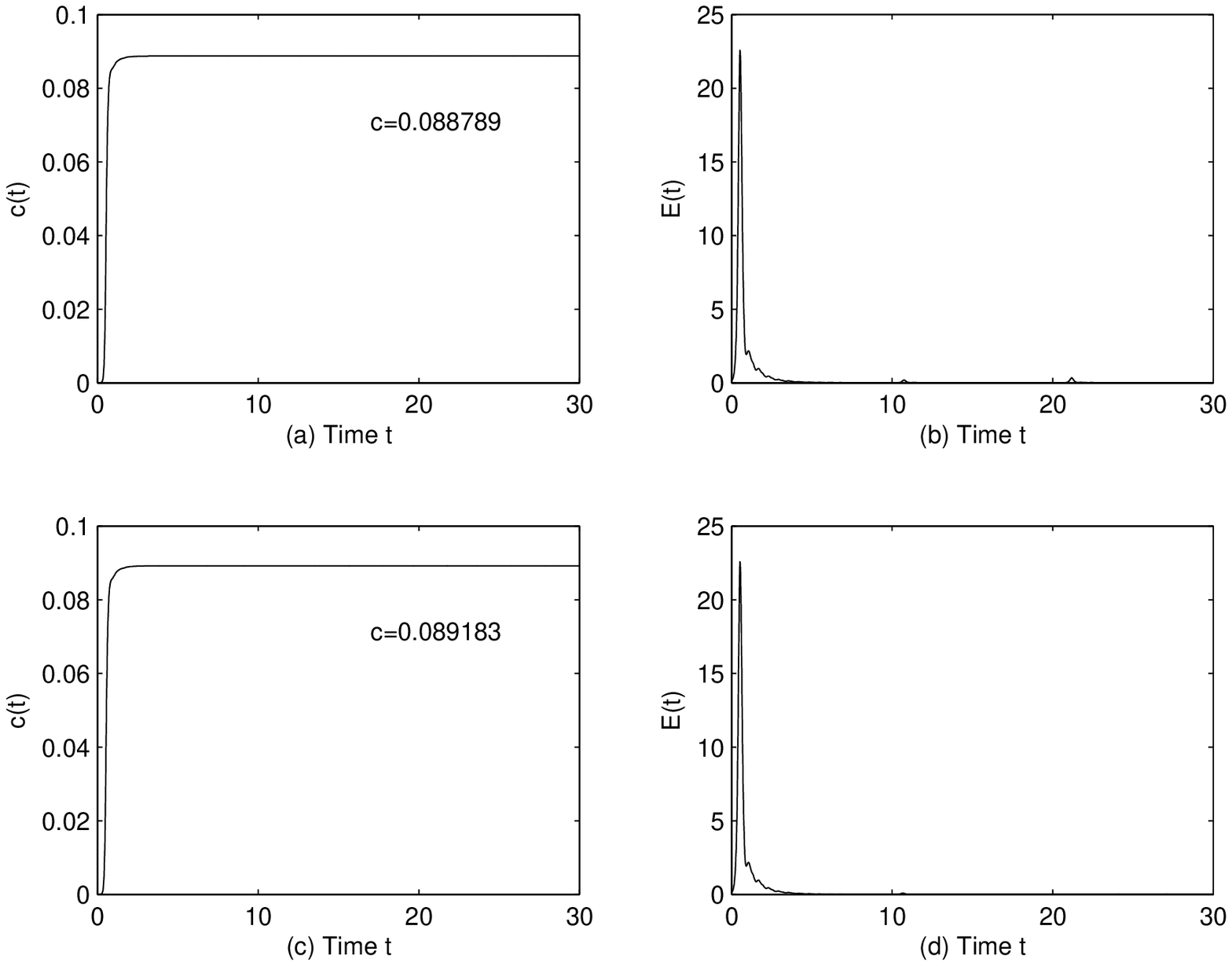}
\end{center}
\begin{center}
Figure 6. The dynamics of $c(t)$ and $E(t)$ for (\ref{any}) with
globally and randomly connected $\tilde{A}$
\end{center}
\subsection{R\"ossler's Oscillator}

The uncoupled equation is :
\begin{eqnarray*}
\left\{
\begin{array}{l}
\dot{x}_1=-x_2-x_3\\
\dot{x}_2=x_1+0.2x_2\\
\dot{x}_3=0.2+x_3(x_1-\mu)
\end{array}
\right.
\end{eqnarray*}
where $\mu=5.7$. Figure 7. (a) and (b) show the dynamics of $c(t)$
and $E(t)$ for linearly coupled systems (\ref{equ1}); while (c) and
(d) show the dynamics of $c(t)$ and $E(t)$ for nonlinearly coupled
systems (\ref{equ2}).
\begin{center}
\includegraphics[height=4.5in, width=5.5in]{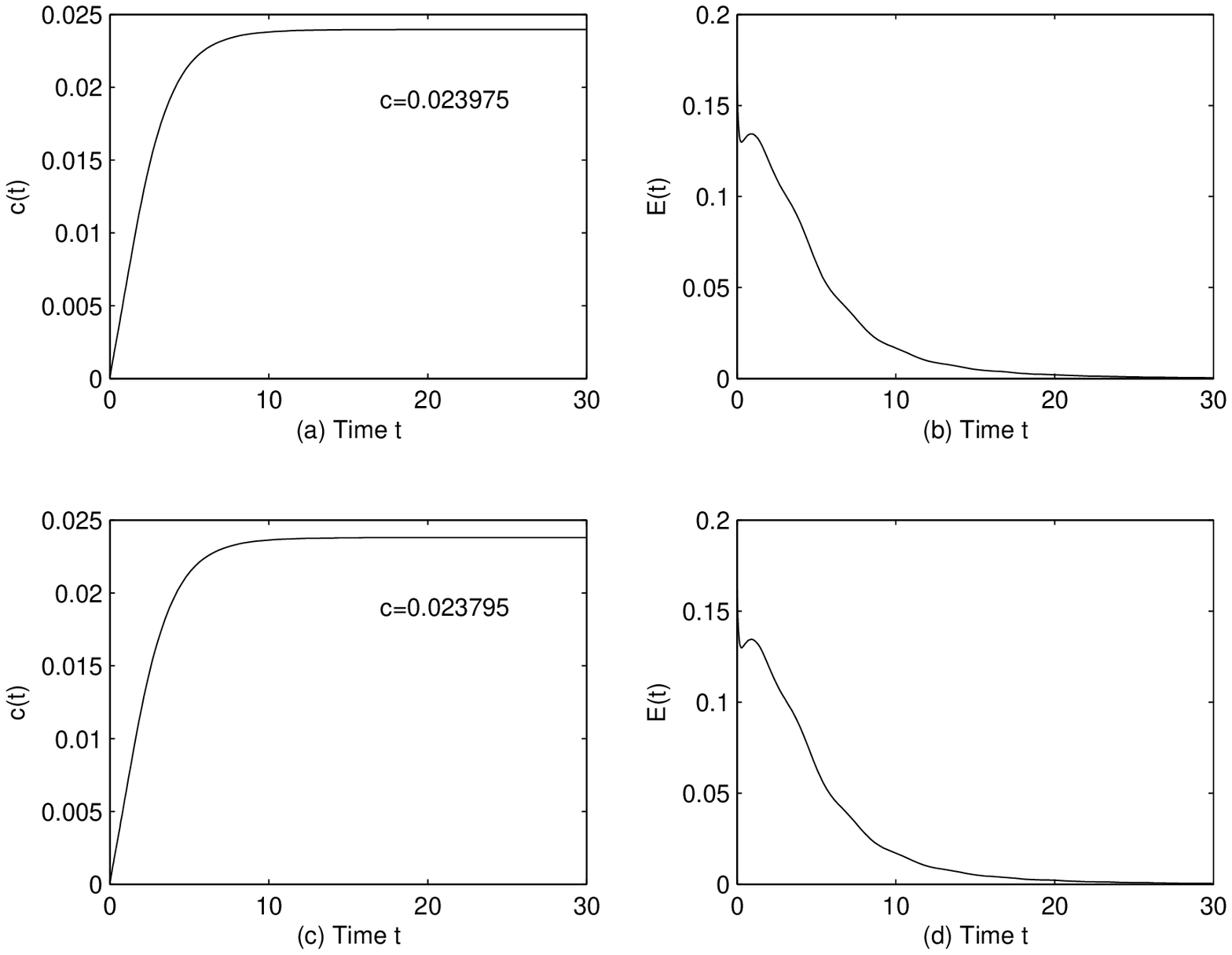}
\end{center}
\begin{center}
Figure 7. The dynamics of $c(t)$ and $E(t)$ for 100 linearly and
nonlinearly\\
coupled R\"ossler's chaotic oscillators with an adaptive coupling
strength
\end{center}

Figure 8. (a) and (b) show the dynamics of $c(t)$ and $E(t)$ for
linearly coupling systems (\ref{any}) when $\tilde{A}$ is the
globally coupled matrix; while (c) and (d) show the dynamics of
$c(t)$ and $E(t)$ for linearly coupled systems (\ref{any}) when
$\tilde{A}\in {\bf A2}$ is a random matrix.
\begin{center}
\includegraphics[height=4.5in,width=5.5in]{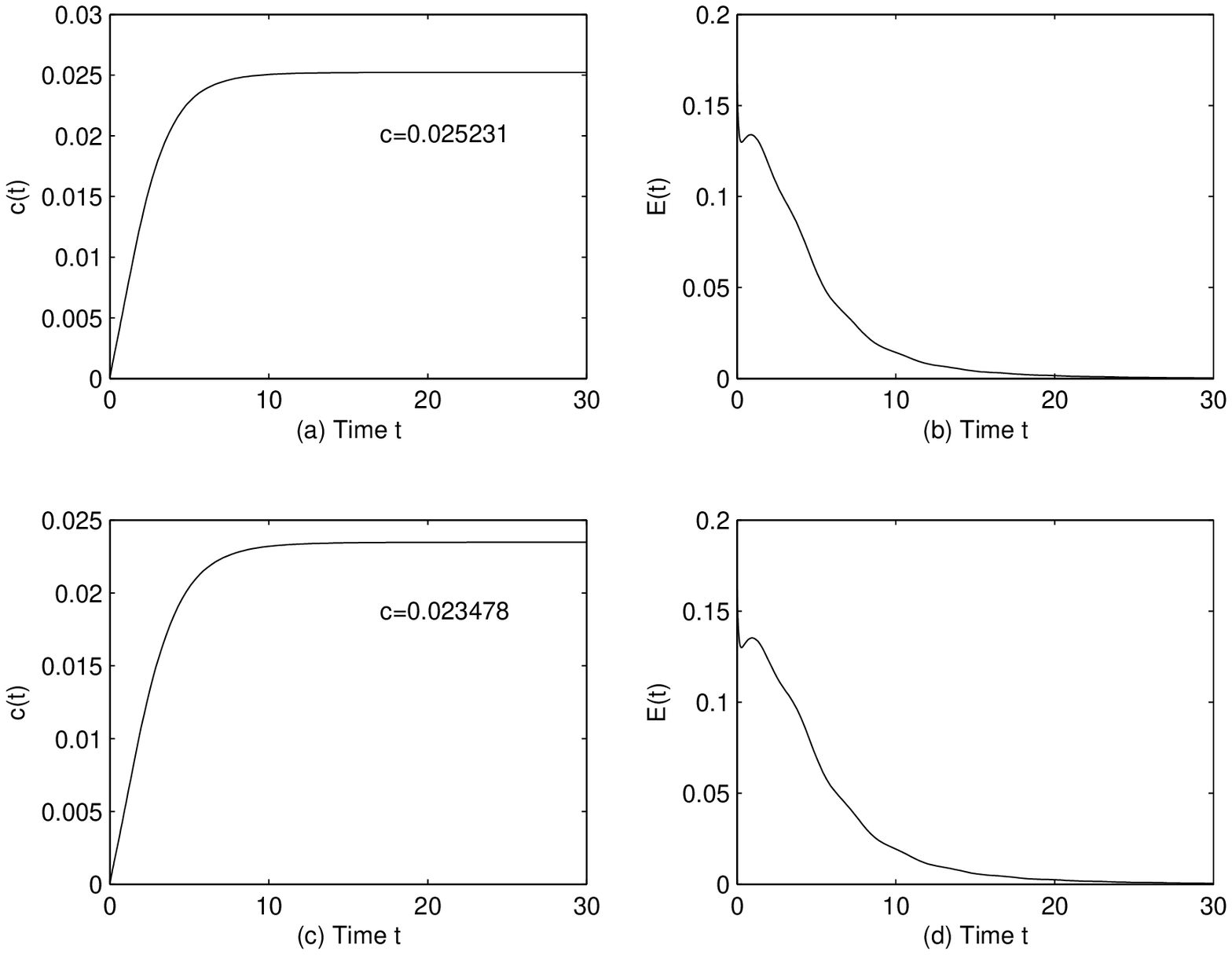}
\end{center}
\begin{center}
Figure 8. The dynamics of $c(t)$ and $E(t)$ for (\ref{any}) with
globally and randomly connected $\tilde{A}$
\end{center}

\section{Numerical Example 2}\quad
In this section, we give some simulations for coupled systems with a
time-varying coupling matrix. As a special case, it includes
node-balanced coupled systems.

We pick
\begin{eqnarray}
A(t)=\left(
\begin{array}{ccc}       p_1       &           &      \\
                                   &p_2        &      \\
                                   &           &p_3
\end{array}\right)
\left(
\begin{array}{ccc}        -5-\sin t-\cos t    &3+\sin t               &2+\cos t         \\
                             2+\cos t         &-5-\sin t-\cos t       &3+\sin t         \\
                             3+\sin t         &2+\cos t               &-5-\sin t-\cos t
\end{array}\right)\label{cmnb}
\end{eqnarray}
as the time-varying coupling matrix, where $p_1,p_2,p_3$ are
positive constants.

It is clear that for all $t$, $A(t)$ has the same left eigenvector
$\xi=(\frac{1}{3 p_1},\frac{1}{3p_2},\frac{1}{3 p_3})^T$
corresponding to eigenvalue 0. In particular, if $p_1=p_2=p_3$, then
the coupling matrix $A(t)$ is node-balanced  for each $t$.

Moreover, it is easy to check that
\begin{eqnarray*}
\lambda_2(t)=-(5+\sin t+\cos t)\leq-(5-\sqrt{2})<0
\end{eqnarray*}
is the largest non-zero eigenvalue of the matrix
\begin{eqnarray*}
\Xi A(t)+A(t)^T\Xi=\frac{5+\sin t+\cos t}{3}\left(
\begin{array}{ccc}  -2 & 1& 1 \\
                     1 &-2& 1 \\
                     1 & 1&-2
\end{array}\right)
\end{eqnarray*}

We couple three Chua's circuits by
\begin{eqnarray}
\left\{\begin{array}{cc}\dot{X}(t)=F(X(t))+c(t){\bf A(t)}X(t),\\
\dot{c}(t)=-\frac{\alpha}{2}X^T(t){\bf \Xi}{\bf A(t)}X(t)
\end{array}\right.\label{time-varying1}
\end{eqnarray}
We also use
$E(t)=\sqrt{(\|x_2(t)-x_1(t)\|^2+\|x_3(t)-x_1(t)\|^2)/2}$ to measure
the error of synchronization.

Figure 9. shows the dynamics of the coupling strength $c(t)$ and
error for different $p_1, p_2, p_3$, Where sub-figures (a) and (b)
denote the dynamics of $c(t)$ and $E(t)$ for $p_1=p_2=1$, $p_3=2$;
and sub-figures (c) and (d) denote the dynamics of $c(t)$ and $E(t)$
for $p_1=p_2=p_3=1$.
\begin{center}
\includegraphics[height=4.6in, width=5.5in]{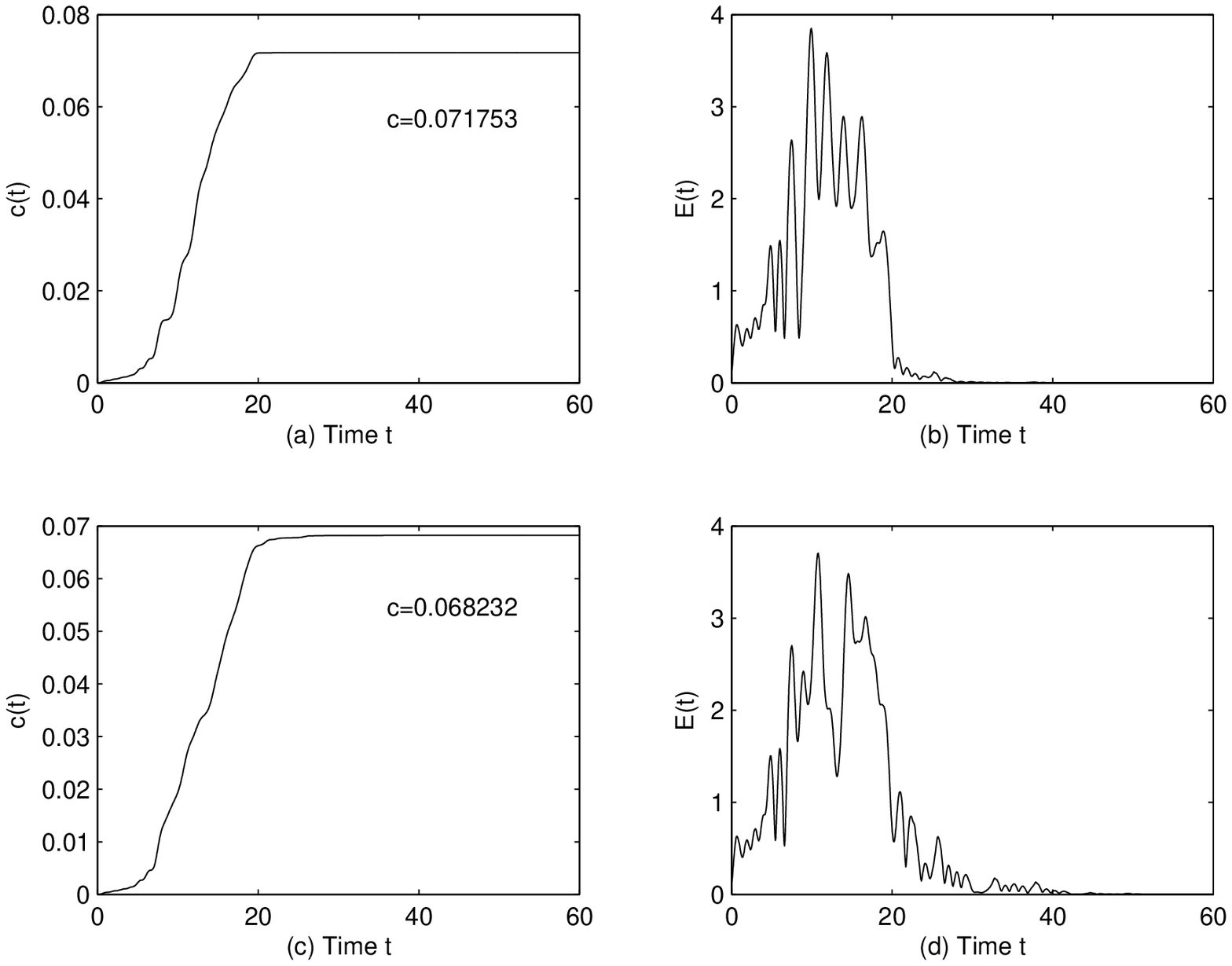}
\end{center}
\begin{center}
Figure 9. The dynamics of $c(t)$ and $E(t)$ for the coupled system
(\ref{time-varying1}).
\end{center}

\section{Conclusions}\quad
In this paper, we propose new algorithms to synchronize linearly or
nonlinearly coupled systems with an adaptive coupling strength.
Unlike those adaptive algorithms existing in the literature, where
coupled systems are synchronized to a special trajectory $s(t)$ or
an equilibrium of the uncoupled system by adding a negative feedback
controller; in this paper, we synchronize linearly and nonlinearly
coupled systems with an adaptive coupling strength without knowing
the synchronization trajectory. By adapting the coupling strength,
we reveal that a large scale of chaotic oscillators can be
synchronized even with a very small coupling strength. It indicates
that chaotic oscillators are very easy to be synchronized.

\section{Acknowledgement}\quad
This project is supported by NSF of China under Grant 60374018,
60574044.

\noindent{\large\bf References}

\begin{description}

\item
Strogatz,S. H. [2001]~ ``Exploring complex networks,'' Nature 410,
pp. 268-276.

\item
Albert R. and Barabasi,  A. L.  [2002]~``Statistical mechanics of
complex networks,'' Rev. Mod. Phys 74, 47-97.

\item
Newman, M. E. J. [2003]~ ``The structure and function of complex
networks,'' SIAM Review 45, 167-256.

\item
Watts D. J. and Strogatz,S. H.  [1998]~ ``Collective dynamics of
small-world,'' Nature 393, 440-442.

\item  Wang, X. F.and Chen, G.  [2002a]~  ``Synchronization in
scale-free dynamical networks: robustness and fragility'', IEEE
Trans. Circuits Syst.-1, 49(1), pp. 54-62.

\item Wang, X. F.and Chen, G. [2002b]~ ``Synchronization
small-world dynamical networks'', Int. J. Bifur. Chaos, 12(1), pp.
187-192.

\item
Mirollo R. E. and Strogatz,  S. H. [1990]~`` Synchronization of
pulse-coupled biological oscillators,'' SIAM J. Appl. Math 50(6),
pp. 1645-1662.

\item
Wei G. W. and Jia,  Y. Q.  [2002]~``Synchronization-based image edge
detection,'' Europhys. Lett 59(6), pp. 814-819.

\item
Xie, Q. X. Chen G. R. and Bollt, E. M. [2002]~ ``Hybrid chaos
synchronization and its application in information processing,''
Math. Comput. Model 35(1-2), pp. 145-163.

\item
Yang, T. and Chua, L. O.  [1997]~``Impulsive control and
synchronization of nonlinear dynamical systems and application to
secure communication,'' Int. J. Bifur. Chaos 7(3), pp. 645-664.

\item
Boccaletti, S. Farini, A.  and Arecchi,  F. T.  [1997]~``Adaptive
synchronization of chaos for secure communication,'' Phys. Rev. E
55(5), pp. 4979-4981.

\item
Liao T. L. and Tsai, S. H. [2000]~``Adaptive synchronization of
chaotic systems and its application to secure communications,''
Chaos, Solitons and Fractals 11, pp. 1387-1396.

\item
Lu, W. L. and Chen,  T. P. [2004]~ ``Synchronization of coupled
connected neural networks with delays,'' IEEE Trans. Circuits Syst.
-I 51, 2491-2503.

\item
Wu and C. W. Chua, L. O. [1995]~``Synchronization in an array of
linearly coupled dynamical systems, '' IEEE Trans. Circuits Syst. -I
42(8), pp. 430-447.

\item
Pecora L. M. and Carroll, T. L. [1998]~ ``Master stability functions
for synchronized coupled systems,'' Phys. Rev. Lett 80(10), pp.
2109.

\item
Belykh, V. N. Belykh  I. V. and Hasler,  M. [2004]~ ``Connection
graph stability method for synchronized coupled chaotic systems,''
Phys. D 195, pp. 159-187.

\item
Belykh, I. Belykh V. and Hasler, M. [2006]~ ``Synchronization in
asymmetrically coupled networks with node balance,'' Chaos 16,
015102.

\item
Wu, C. W. [2005]~ ``Synchronization in networks of nonlinear
dynamical systems coupled via a directed graph,'' Nonlinearity 18,
pp. 1057-1064.

\item
Lu W. L. and Chen,  T. P. [2006]~ ``New approach to synchronization
analysis of linearly coupled ordinary differential systems,'' Phys.
D 213, pp. 214-230.

\item
Pikovsky, A. Rosenblum, M. and Kurths,  J. [2001]~ Synchronization:
A Universal Concept in Nonlinear Sciences, Cambridge University
Press.

\item
Chen T. P. and Zhu,  Z. M. [2007] ~``Exponential synchronization of
nonlinear coupled dynamical networks,'' Int. J. Bifur. Chaos (to
appear).

\item
Ott, E. Grebogi C. and Yorke,  J. A.  [1990]~ ``Controlling chaos
,'' Phys. Rev. Lett 64 (11), pp. 1196-1199.

\item
Chen M. Y. and Zhou,  D. H. [2006]~ ``Synchronization in uncertain
complex networks,'' Chaos 16, 013101.

\item
Zhou, J. Lu  J. N. and L\"{u}  J. H. [2006]~ ``Adaptive
synchronization of an uncertain complex dynamical network,'' IEEE
Trans. Automatic Contr  51(4), pp. 652-656.

\item
Wu, C. W. Yang T. and Chua, L. O.  [1996]~ ``On adaptive
synchronization and control of nonlinear dynamical systems,'' Int.
J. Bifur. Chaos 6(3), pp. 455-471.

\item
Dedieu H. and Ogorzalek,  M. J. [1997]~``Identifiability and
identification of chaotic systems based on adaptive
synchronization,'' IEEE Trans. Circuits Syst. -I 44(10), pp.
948-962.

\item
Hong, Y. G. Qin, H. S. Chen,  G. R. [2001]~``Adaptive
synchronization of chaotic systems via state or output feedback
control,'' Int. J. Bifur. Chaos 11(4), pp. 1149-1158.

\item
Lian, K. Y. Liu,  P. Chiang T. S. and Chiu, C. S.  [2002]~``Adaptive
synchronization design for chaotic systems via a scalar driving
signal,'' IEEE Trans. Circuits Syst. -I 49 (1), pp. 17-27.

\item
Li Z. and Chen,  G. R.  [2004]~``Robust adaptive synchronization of
uncertain dynamical networks,'' Physics Letters A  324, pp. 166-178.

\item
Elabbasy, E. M. Agiza  H. N. and El-Dessoky,  M. M. [2005]~``Global
synchronization criterion and adaptive synchronization for new
chaotic system,'' Chaos, Solitons and Fractals 23, pp. 1299-1309.

\item
Yau, H. T. Lin  J. S. and YAN,  J. J.  [2005]~``Synchronization
control for a class of chaotic systems with uncertainties,''Int. J.
Bifur. Chaos 15 (7), pp. 2235-2246.

\item
Cao J. D. and Liu,  J. Q. [2006]~ ``Adaptive synchronization of
neural networks with or without time-varying delay,'' Chaos 16,
013133.

\item
Zhou, J. Chen T. P. and Xiang,  L.  [2006]~``Robust synchronization
of delayed neural networks based on adaptive control and parameters
identification,'' Chaos, Solition and Fractals 27, pp. 905-913.

\item
Zhang, H. G. Huang, W.  Wang  Z. L. and Chai,  T. Y.
[2006]~``Adaptive synchronization between two different chaotic
systems with unknown parameters,'' Physics Letters A 350, pp.
363-366.

\item
S. Boccaletti, V. Latora, Y. Moreno, M. Chavez and D.-U. Hwang,
[2006]~ ``Complex networks: Structure and dynamics,'' Physics
Reports 424, pp. 175-308.

\item
Kuang, Y. [1993]~"Delay Differential Equations in Populaiton
Dynamics," New York: Academic Press.

\end{description}

\end{document}